\documentclass[amscd]{amsart}

\usepackage{amscd,amssymb,amsmath,amstext}

\newtheorem{Thm}{Theorem}
\newtheorem{Cor}{Corollary}
\newtheorem{Lem}{Lemma}
\newtheorem{Prop}{Proposition}

\newtheorem{Fact}{Fact}

\theoremstyle{remark}
\newtheorem{Rem}{Remark}
\newtheorem{Def}{Definition}

\newcommand{\cal}{\mathcal}

\newcommand{\oplusl}{\bigoplus\limits}
\newcommand{\cupl}{\bigcup\limits}
\newcommand{\suml}{\sum\limits}

\newcommand{\bu}{\bullet}

\newcommand{\const}{{\rm const}}

\newcommand{\gr}{{\bf gr}}

\newcommand{\Lh}{{^L{\mathfrak h}}}

\newcommand{\Lg}{{^L{\mathfrak g} }}
\newcommand{\Lgtil}{{^L \tilde {\mathfrak g} } }

\newcommand{\Gr}{{{\cal G}\frak r}}
\newcommand{\Fl}{{{\cal F}\ell}}

\newcommand{\stil}{{\tilde{s}}}
\renewcommand{\c}{{\underline c}}
\newcommand{\LG}{{ ^LG}}
\newcommand{\LB}{{ ^LB}}                                  
\newcommand{\Xcheck}{{ X\check{\,}}}

\newcommand{\G}{{\cal G}}
\newcommand{\F}{{\cal F}}
\newcommand{\imbed}{\hookrightarrow}
\newcommand{\Ad}{{\rm Ad\,}}
\newcommand{\Aut}{{\rm Aut\,}}
\newcommand{\Fq}{{{\Bbb F}_q}}
\newcommand{\Fqbar}{{\overline{\Bbb F}_q}}
\newcommand{\Ql}{{\overline {{\Bbb Q}_l} }}
\newcommand{\Qlbar}{{\overline {{\Bbb Q}_l} }}

\newcommand{\Zet}{{\Bbb Z}}
\newcommand{\Ce}{{\Bbb C}}
\newcommand{\kk}{ k }
\newcommand{\I}{{{\mathbb I}_\c}}
\newcommand{\II}{{\mathbb I}}

\newcommand{\bG}{{\bf G}}
\newcommand{\bI}{{\bf I}}
\newcommand{\bK}{{\bf K}}

\renewcommand{\P}{{\cal P}}

\newcommand{\To}{\longrightarrow}
\newcommand{\iso}{{\widetilde \longrightarrow}}
\newcommand{\sur}{ \twoheadrightarrow}
\newcommand{\A}{{\cal A}}
\newcommand{\B}{{\cal B}}

\renewcommand{\H}{{\cal H}}

\newcommand{\epf}{\square}
\newcommand{\gal}{{\check{\ }}}
\newcommand{\cO}{{\cal O}}
\renewcommand{\O}{{\cal O}}
\newcommand{\J}{{\cal J}}
\newcommand{\N}{{\cal N}}
\newcommand{\Ntil}{{\tilde {\cal N}}}
\newcommand{\R}{{\cal R}}

\newcommand{\Rk}{{\cal R}^k}
\newcommand{\unO}{{\underline{\cal O}}}
\newcommand{\Jf}{J^f}
\newcommand{\Af}{A^f}
\newcommand{\mon}{{\frak M}}

\newcommand{\Stt}{{\bf St}}

\newcommand{\chitil}{{\tilde \chi}}

\renewcommand{\proof}{{\it Proof }}

\def\square{\hbox{\vrule\vbox{\hrule\phantom{o}\hrule}\vrule}}

\begin{document}

\title[]{On tensor categories attached to cells in affine Weyl groups}
\author{Roman Bezrukavnikov}

\begin{abstract}
This note is devoted to Lusztig's bijection between unipotent
conjugacy classes in a simple complex algebraic group and
 2-sided cells in the affine
Weyl group of the Langlands dual group; and also to
the description of the reductive quotient of the centralizer of the
unipotent element in terms of convolution of perverse sheaves on affine flag
variety of the dual group conjectured by Lusztig in  \cite{conje}.
Our main tool is a recent construction
by Beilinson,  Gaitsgory and Kottwitz,
 the so-called sheaf-theoretic construction of
the
center of an affine Hecke algebra (see \cite{KGB}).
We show how this remarkable construction provides a geometric interpretation of
the bijection, and allows to prove the conjecture.
\end{abstract}
\maketitle

{\bf Acknowledgement.} I am much indebted to Michael Finkelberg;
among the many things he taught me are the theory of cells in Weyl groups and
Lusztig's conjecture (partially proved below).
This note owes a lot to Dennis Gaitsgory, who
 explained to me the ``sheaf-theoretic center'' construction,
and made some helpful suggestions.
 I also thank Alexander Beilinson and Vladimir Drinfeld
for useful comments, and Viktor Ostrik for stimulating interest.
The results of this note were
 obtained (in a preliminary form)
during the author's participation in the special year
on Geometric Methods in Representation Theory (98/99)
at IAS; I thank IAS for
its hospitality, and NSF grant for financial support.

\begin{section}{Introduction.}
Let $G$ be a split
simple algebraic group. Let $W$ be the corresponding  affine
Weyl group; and let $J$ be the
asymptotic affine Hecke algebra \cite{Cells2}.
Recall that  $J$ is an algebra over $\Zet$,
and it comes with a basis $t_w$ parametrzied by $W$.
The group $W$ is the union of two-sided
cells \cite{Cells0}, and the algebra $J$ is the direct sum of algebras,
$J=\oplus J_\c$, where $\c$ runs over the set of two-sided cells,
and $J_\c$ is the span of $t_w$, $w\in \c$.
Let also $W_f\subset W$ be the finite Weyl group,
and $W^f\subset W$ be the set of minimal length representatives
of double cosets $W_f\backslash W/W_f$. Let $\Jf$ be the span of
$t_w$, $w\in W^f$, and $\Jf_\c=\Jf\cap J_c$.
It follows from the result  of \cite{LXi}
that $\Jf\subset J$, $\Jf_\c
\subset J_\c$ are  subalgebras.

Let $\LG$ be the Langlands dual group (over an algebraically
closed field of characteristic zero).
 In \cite{Cells4} Lusztig constructed
 a bijection between two-sided cells
in $W$ and unipotent conjugacy classes in $\LG$; for a 2-sided cell
$\c$ we will  denote  by $N_\c\in \LG$  a representative of the unipotent
conjugacy class corresponding to $\c$.
In \cite{conje}
the based algebras $J_\c$, $\Jf_\c$ are realized as Grothendieck groups
of certain semisimple
monoidal (tensor without commutativity) categories, which we denote
respectively by $A_\c$,
$\Af_\c$  (thus $\Af_\c$ is a monoidal subcategory in $A_\c$).
 Categories $A_\c$, $\Af_\c$ are defined as subcategories of
(semisimple) perverse sheaves on the affine flag manifold of $G$;
 the monoidal structure is provided by the truncated
 convolution, see {\it loc. cit.}

 The following  conjectural description of monoidal
categories $A_\c$, $\Af_\c$ (and hence of based rings $J_\c$, $J_\c^f$)
was proposed
in \cite{conje}, \S 3.2.
 Let $Z=Z_\LG(N_\c)$ be the centralizer of $N_\c$ in $\LG$.
For any two-sided cell $\c$
there exists (conjecturally) a finite set $X$ with a $Z$ action,
such that $A_\c$ is equivalent to the category $Vect^Z_{ss}(X\times X)$
of semisimple $Z$-equivariant
sheaves on $X\times X$, with monoidal structure given by convolution.
The finite set $X$ should contain a preferred point $x\in X$ fixed by $Z$;
thus the category $Vect^Z_{ss}(X\times X)$ contains a monoidal subcategory
$Vect^Z_{ss}(\{(x,x)\})=Rep^{ss}(Z)$ (the category of semisimple
representations of $Z$). This subcategory should be identified with $\Af_\c
\subset A_\c$.

In this note we (extend and) prove the part of the above
conjecture, which asserts the equivalence $\Af_\c=Rep^{ss}(Z)$.
The proof is based on a recent construction by Gaitsgory (following an idea
of Beilinson and
 Kottwitz), see \cite{KGB}.

Recall that the link between representations of $\LG$ and theory of perverse
sheaves on affine Grassmanian  of $G$ is provided by the so-called geometric
version of the
Satake isomorphism (the idea going back to \cite{doistor} is developed
in \cite{Gi} and \cite{MV}; see also \cite{BD}). The classical Satake
isomorphism is an isomorphism between the spherical Hecke
algebra $\H_{sph}$ and the ring $\R(\LG)$, where for an algebraic
group $H$ we write $\R(H)$ for its representation ring (Grothendieck
group of the category $Rep(H)$ with multiplication provided by the tensor
product). Its geometric (or categorical)
version is an  equivalence of
tensor categories between
 $Rep(\LG)$ and perverse sheaves on the affine Grassmanian.

 Further, a theorem of Bernstein (see  e.g. \cite{doistor}, Proposition
8.6) asserts that
the center of the Iwahori-Matsumoto Hecke algebra $\H$ is also
isomorphic to  $\R (\LG)$.
 The main result of \cite{KGB} provides a geometric (categorical)
counterpart of this  isomorphism. More precisely, it defines an action of the
tensor category $Rep(\LG)$ on the category of perverse sheaves on affine flags;
 on the level of Grothendieck groups this
 amounts to the action of the center of $\H$, $Z(\H)=\R(\LG)$
on $\H$.
This action enjoys various favorable properties; it also carries a canonical
unipotent automorphism $\mon$ (the monodromy).

The idea of the present note is that one can identify certain
subquotient categories  of perverse sheaves on affine flags
 related to a 2-sided cell $\c$ with
$Rep(H)$ for a subgroup $H\subset Z(N_\c)$, in such a way that the canonical
action of $Rep(\LG)$ constructed in \cite{KGB} is identified with the
tautological action of $Rep(\LG)$ on $Rep(H)$
$$V:W\mapsto Res^\LG_H (V)\otimes W;$$
 moreover this requirement fixes
the identification with $Rep(H)$ uniquely.
 Then the monodromy automorphism  $\mon$ provides (by Tannakian formalism)
a unipotent element $N\in \LG$, commuting with $H$. This element lies in the
conjugacy class attached to $\c$ by Lusztig.

With these tools in hand, the Lusztig's conjecture becomes an
exercise in Tannakian formalism (at least
modulo some powerful Theorems of
Lusztig on the structure of asymptotic Hecke algebras).

I want to point out that although this note contains  a complete
proof of the stated result, the argument presented here may not be optimal.
In particular, at several places we use ``as a black box''
results of Lusztig on asymptotic Hecke algebras to check a categorical
property of perverse sheaves on affine flags.
In Remarks \ref{op1}, \ref{op2} below
   we discuss a possible plan for replacing some of these
uses by a direct geometric argument (in other words, for providing a geometric
proof of Lusztig's results). However, we do not suggest such a
geometric proof e.g. for  Lusztig's description of the unit element in the
 asymptotic Hecke algebra (equivalent to certain properties of
 Duflo involutions).

\section{Preliminaries on tensor categories}
In this section we collect some general Lemmas on Tannakian formalism.
They are probably well-known (anyway, the proof is straightforward)
but a reference was not found.

For a category $A$ (full subcategory of a triangulated category) 
 we will
write $K(A)$ for the Grothendieck group of $A$. For an object $X\in A$
we will denote its class by $[X]\in K(A)$.

\begin{Lem}\label{tochno}
 Let $\A$ be an abelian             
category with  all objects  having
finite length. Let $\otimes$ denote a functor $\A\times \A \to A$ which
is linear and mid-exact in each variable.
For an object $X\in \A$ let $X^{ss}$ denote
the semisimplification of $X$. Let $X, Y\in \A$ be two objects satisfying
\begin{equation}\label{eqtochno}
[X\otimes Y]=[X^{ss}\otimes Y^{ss}]
\end{equation}
where $[X]$ denotes the class of $X$ in the Grothendieck group  $K^0(\A)$.
Then \eqref{eqtochno} remains true when $X,Y$ are replaced by a
subquotient. Moreover,
 if $X'$, $Y'$ are subquotients of respectively $X$, $Y$
and  $0\to X'' \to X'\to X'''\to 0$ is an exact sequence, then
the sequence $0\to X''\otimes Y'\to X'\otimes Y'\to X'''\otimes Y'\to 0$
is exact.
\end{Lem}
\proof
For $\alpha,\beta \in K^0(\A)$ let us write $\alpha \leq \beta$
if $\beta-\alpha$ is a class of actual  (as opposed to virtual)
object. Then for any s.e.s. $ 0\to X_1 \to X_2\to X_3\to 0$
we have
\begin{equation}\label{nerav}
[X_2\otimes Y]\leq [X_1\otimes Y]+[X_3\otimes Y]
\end{equation} for any $Y$,
with equality being true iff the sequence   $ 0\to X_1\otimes Y \to
X_2\otimes Y \to X_3\otimes Y\to 0$ is exact. Hence the first statement
implies the second.

By induction in Jordan-Hoelder series we see that $[X\otimes
Y] \leq [X^{ss}\otimes Y^{ss}]$ for any $X,Y$. Moreover, if  for
subquotients $X',Y'$ of $X,Y$   the strict
inequality    $[X'\otimes
Y'] < [{X'}^{ss}\otimes {Y'}^{ss}]$ holds, then successive use of \eqref{nerav}
shows that $[X\otimes Y]< [X^{ss}\otimes Y^{ss}]$, which contradicts
\eqref{eqtochno}.
\epf

\begin{Def}\label{defcen}
Let $\A$ be a monoidal category,
 and $\B$ be a tensor (symmetric monoidal)
category. A {\it central functor} from $\B$ to $\A$
 is a monoidal functor $F:\B\to \A$ together with an isomorphism
\begin{equation}\label{abscen}
\sigma_{X,Y}: \ F(X)\otimes Y \cong Y\otimes F(X)
\end{equation}
fixed for all $X\in \B$, $Y\in \A$, subject to the following compatibilities.

i) $\sigma_{X,Y}$ is functorial in $X,Y$;

ii) For $X,X'\in \B$ the isomorphism $\sigma _{X,F(X')}$ coincides with
the composition
$$F(X)\otimes F(X')\cong F(X\otimes X')\cong F(X'\otimes X)
\cong F(X')\otimes F(X)$$
(where the middle isomorphism comes from the commutativity constraint
in $\B$, and the other two from the tensor structure on $F$).

iii)  For $Y_1, Y_2\in \A$ and $X\in \B$
 the composition
$$\begin{CD}
F(X)\otimes Y_1\otimes Y_2
 @>{\sigma_{X,Y_1}\otimes Y_2}>>
Y_1 \otimes F(X) \otimes Y_2 @>{Y_1\otimes \sigma_{X,Y_2}}>>
Y_1\otimes Y_2\otimes F(X)
\end{CD}
$$
coincides with $\sigma_{X, Y_1\otimes Y_2}$.

iv)   For $Y\in \A$ and $X_1,X_2\in \B$
the composition
$$
\begin{CD}
F(X_1\otimes X_2)\otimes Y\cong F(X_1)\otimes F(X_2)\otimes Y
@>{F(X_1)\otimes    \sigma_{X_2,Y}}>>
F(X_1) \otimes Y \otimes F(X_2)
@>{ \sigma_{X_1,Y}\otimes F(X_2)}>> \\
Y\otimes F(X_1)\otimes F(X_2)=Y\otimes F(X_1\otimes X_2)
\end{CD}
$$
coincides with $\sigma_{X_1\otimes X_2, Y}$.

\end{Def}

\begin{Rem}
\footnote{The content of this remark was communicated to me by
 Drinfeld.}
For a monoidal category $\A$ its {\it center} ${\mathcal Z}(\A)$
is defined as the category
of pairs $(X, \sigma)$, where $X\in \A$, and $\sigma=(\sigma_Y)$
is a collection of isomorphisms
$\sigma_Y: X \otimes Y  \cong Y\otimes X$, subject to certain compatibilities.
Then ${\mathcal Z}(\A)$
 turns out to carry a natural structure of a braided monoidal
category. If $\A$ is the category of representations of a Hopf algebra $A$,
 then ${\mathcal Z}(\A)$ is identified with the category of representations
of Drinfeld's double of $A$.

One can check  
that a central functor from a (symmetric) tensor category
$\B$ to a monoidal category  $\A$ is the same
as a monoidal functor from $\B$ to ${\mathcal Z}(\A)$, which intertwines
commutativity in $\B$ with the braiding in  ${\mathcal Z}(\A)$.
\end{Rem}

\begin{Prop}\label{genlem} Let $k$ be an algebraically closed field.
 Let $\A$ be a $k$-linear abelian
monoidal category  with a unit object $\II$ such that $End(\II)=k$.
We assume that the product in $\A$ is exact in each variable.

Let $G$ be an algebraic group over $k$, and $Rep(G)$ be the category
of its finite dimensional algebraic representations.
Let $F:Rep(G)\to \A$ be an exact central functor.
Suppose that any object $Y\in \A$ is isomorphic to a subquotient of
$F(X)$ for some
 $X\in Rep(G)$.

Assume also\footnote{These assumptions are not  necessary, and are
imposed to shorten the proof.}
 that $k$ is uncountable and $Hom(X,Y)$ is finite dimensional for
$X,Y\in \A$.
 Then there exists an algebraic subgroup $H\subset G$,
and an equivalence of monoidal categories $\Phi: Rep(H)\iso \A$,
 such that
$$F\cong \Phi \circ Res^G_H.$$

The subgroup $H\subset G$ is defined uniquely up to conjugation.
\end{Prop}

\proof
$G$ acts on itself by left translations, making
the space $\O(G)$ of regular functions on $G$ an algebraic $G$-module, thus an
 ind-object of $Rep(G)$. Let $\unO _G$ denote this ind-object.
It is a ring ind-object, i.e. we have a multiplication morphism
$m:\unO_G \otimes \unO_G \to \unO_G$ satisfying the usual commutative
ring axioms. For a ring ind-object $?$ we will write $m_?$ for the
multiplication morphism $m_?:?\otimes ? \to ?$.
Let $\J  \subsetneq
 F(\unO_G)$ be a maximal left ideal subobject, i.e. $\J$ is a maximal proper
 ind-subobject in $F(\unO_G)$
 satisfying
\begin{equation}\label{ideal_1}
m_{F(\unO_G)}(F(\unO_G) \otimes \J)\subset \J.
\end{equation}
Then $\J$ is also a right
ideal subobject, i.e. we have
\begin{equation}\label{ideal_2}
m(\J\otimes F(\unO_G) )\subset \J;
\end{equation}
indeed, commutativity of $\unO_G$, and property (ii) in
the definition of a central functor
show that $$m_{F(\unO_G)} \circ \sigma_{\unO_G, F(\unO_G)} =
m_{F(\unO_G)},$$
and property (i) yields the equality
$$m_{F(\unO_G)}|_{\J\otimes F(\unO_G)} \circ \sigma_{\unO_G, \J} =
m_{F(\unO_G)}|_{F(\unO_G) \otimes \J},$$
which implies \eqref{ideal_2}.

Set $\unO_H=F(\unO_G)/\J$. \eqref{ideal_1}, \eqref{ideal_2} imply that
$\unO_H$ is a ring ind-object of $\A$.
Thus the category of $\unO_H$-module (ind)objects in $\A$
 is well-defined. We will denote this category by $\unO_H-mod$,
 call its objects
$\unO_H$-modules, and write
$Hom_{\unO_H}$ instead of $Hom_{\unO_H-mod}$.

Then $\unO_H-mod$ is an abelian category, and $\unO_H\in \unO_H-mod$
is a simple object. Thus $K=End_{\unO_H}(\unO_H)$ is a division algebra, and
$V\mapsto V\otimes _K \unO_H$ is an equivalence between (right)
finite $K$-modules
and the full subcategory in $\unO_H-mod$ generated by $\unO_H$ under
finite
direct sums and subquotients (we will call such $\unO_H$-modules free;
argument below implies that in fact any $\unO_H$ module is free).

\begin{Lem}\label{Kk}
We have $K=Hom_\A(\II, \unO_H)=k$.
\end{Lem}

\proof  $\unO_H$ is a  unital ring object since $\unO_G$ is, i.e.
the unit $\iota:\II \imbed \unO_H$ is fixed.
The map $\phi\mapsto \phi\circ \iota$ provides an isomorphism
$End_{\unO_H}(\unO_H)\cong Hom_\A (\II,\unO_H)$, with inverse isomorphism
given by $i\mapsto m_{\unO_H}
\circ (i\otimes Id_{\unO_H})$. Thus the first equality
is clear.

To  check the second one, notice that cardinality of a basis of
a division algebra $K$ over an algebraically closed field $k$ is not less
than the cardinality of $k$ (indeed, for $x\in K$, $x\not \in k$ the elements
$(x-\lambda)^{-1}$, $\lambda \in k$ are linearly independent).
However, $\unO_G$ is a countable union of objects of $Rep(G)$, hence
$Hom_\A(\II, F(\unO_G))$ is at most countable dimensional. \epf

\begin{Cor}\label{fiba}
 a) For any $X\in \A$ we have $X\otimes \unO_H\cong V\otimes
\unO_H$ for some finite dimensional $k$-vector space $V$.

b) The functor $\Phi_H :\A \to Vect$
given by $\Phi_H: X\mapsto Hom(\II, X\otimes \unO_H)$
is exact, admits a structure of a monoidal functor, and we have a canonical
isomorphism of monoidal functors
\begin{equation}\label{PhiGH}
\Phi_G \cong \Phi_H \circ F,
\end{equation}
where $\Phi_G:Rep(G)\to Vect$ is the fiber functor.
\end{Cor}

\proof a) Let first $X=F(Y)$ for $Y\in Rep(G)$.
  We have an isomorphism of $\unO_G$-modules
 $Y\otimes \unO_G\cong \Phi_G(Y)\otimes \unO_G$, hence an isomorphism
of $F(\unO_G)$-modules
$$X\otimes F(\unO_G) \cong \Phi_G(Y)\otimes F(\unO_G).$$
Replacing each side of the last equality by
 the maximal quotient on which $F(\unO_G)$ acts through $\unO_H$ we get
an isomorphism of $\unO_H$-modules:
$$F(Y)\otimes \unO_H \cong \Phi_G(Y)\otimes \unO_H.$$
Since any $X\in \A$ is a subquotient of $F(Y)$ for some $Y\in Rep(G)$
we see that the $\unO_H$-module
$X\otimes \unO_H\in \unO_H-mod$ is a subquotient of the free $\unO_H$-module
$\Phi_G(Y)\otimes \unO_H$, hence is also free (i.e. has the form $ V\otimes
\unO_H$ for a $k$-vector space $V$).

Proof of (b). Notice that (a) together
with Lemma \ref{Kk} imply that $X\otimes \unO_H = \Phi_H(X)
\otimes \unO_H$ canonically for  $X\in \A$. This shows exactness of $\Phi_H$,
and also establishes monoidal structure on $\Phi_H$, for we have
$$X\otimes Y\otimes \unO_H=X\otimes (\Phi_H(Y)\otimes \unO_H)=
\Phi_H(X)\otimes \Phi_H(Y) \otimes \unO_H.$$
Finally, the isomorphism $$X\otimes \unO_G\cong \Phi_G(X)\otimes \unO_G$$
for $X\in Rep(G)$ yields (by applying $F$, and taking the maximal quotient
on which $\unO_G$ acts through $\unO_H$) an isomorphism $$F(X)\otimes \unO_H
\cong \Phi_G(X)\otimes \unO_H,$$
hence  an isomorphism \eqref{PhiGH}.
\epf

We can now finish the proof of Proposition \ref{genlem}. According to
\S 2 of \cite{DM}, a functor $\Phi_H$ as in
Corollary \ref{fiba} (b) above
yields  a bialgebra $A$, an equivalence
of monoidal categories $\Psi_H: Comod_A \cong \A$, a morphism of bialgebras
$\phi: \O(G) \to A$ and an isomorphism of monoidal functors
$$F\cong  \Psi_H \circ \phi_*.$$
Since any object of $\A$ is a subquotient of $F(X)$ for some $X\in Rep(G)$,
the morphism $\phi:\O(G) \to A$ is surjective. Thus $A=\O(H)$ for a Zarisski
closed subsemigroup $H\subset G$. Thus Proposition \ref{genlem}
follows from the next  Lemma. \epf

We formulate the Lemma in a slightly greater generality than needed for our
application (since the proof is the same).

\begin{Lem} Let $G$ be a  group scheme of finite type over a Noetherian
ring $k$. Then
a closed subsemigroup scheme $H\subset G$
 is  a subgroup scheme.
\end{Lem}

\proof\footnote{I thank Dima Arinkin, to whom this proof is due.}
 We should check that for any commutative $k$-algebra
 $R$ the subsemigroup
$H(R)\subset G(R)$ is a subgroup. It is enough to check this for
$R$ of finite type over $k$. For $g\in H(R)$ let $L_g:G_R\to G_R$
be the (left)  multiplication by $g$ (here
the subindex $_R$ denotes  base change to $R$).
 Then $L_g(H_R)\subseteq H_R$; we have
to check that in fact $L_g(H_R)=H_R$. But otherwise $H_R\supsetneq L_g(H_R)
\supsetneq L_g^2(H_R) \supsetneq \dots$ is an infinite decreasing chain of
closed subschemes in $G_R$, which constradicts the fact that $G_R$ is
Noetherian. \epf

\begin{Rem}
In this Remark we outline an alternative argument, which is shorter
than the proof of Proposition \ref{genlem} presented above, but
uses a deep Theorem of Deligne \cite{De}, and proves a weaker statement.

In the situation of Proposition  \ref{genlem}
assume that $char(k)=0$, and also that rigidity on the target category
$\A$ is given.\footnote{We do not know whether it is true that
for an abelian tensor category
with exact tensor product a subquotient of a rigid object of finite length
 is rigid.
If the answer to this question is positive, we can drop here the assumption
 that rigidity on $\A$ is given.}
 (In view of Remark \ref{rig_c} below
 this weaker statement is sufficient
for our application.) Then one can show first that
there exists a unique commutativity constraint on $\A$ compatible with one
in $\B$; thus $\B$ is a Tannakian category. Now a Theorem of Deligne
(\cite{De} Theorem 7.1) says that for an algebraically closed
 field $k$ of characteristic 0,
a $k$-linear Tannakian category $\A$ admits a fiber functor,
 and is identified with the category of representations of an algebraic group,
provided for any object $X\in \A$ we have $\Lambda^n(X)=0$
for large $n$ (where $\Lambda^n$ stands the $n$-th exterior power).
If an object $Y$ is a subquotient of $X$, then $\Lambda^n(Y)$ is a subquotient
of $\Lambda^n(X)$, in particular $\Lambda^n(X)=0 \Rightarrow \Lambda^n(Y)=0$.
Thus $\A$ satisfies conditions of Deligne's Theorem, and hence $\A=Rep(H)$
for an algebraic group $H$ by that Theorem. The tensor
functor $F:Rep(G)\to \A= Rep(H)$ yields by Tannakian formalism a homomorphism
$H\to G$; since any object of $Rep(H)$ is a subquotient of $F(X)$ for
$X\in Rep(G)$, this homomorphism is injective.
\end{Rem}

\end{section}

\begin{section}{Notations and recollections}
\subsection{General notations}
Let $G$ be a split
simple algebraic group over $\Zet$.
Let $F=\Fq((t))$ be a local field of
 prime characteristic, and
$O=\Fq[[t]]$ be its ring of integers.\footnote{As usual one
could replace $\Fq$ by an algebraically closed field of characteristic
 zero; then we would have to work with Hodge $D$-modules instead of Weil
sheaves in the proof of Lemma \ref{phiBV} below.}
 Let $I\subset G(O) \subset
G(F)$ be the Iwahori subgroup.
 There exist
canonically defined group schemes $\bK$, $\bI$ over $\Fq$ (of infinite type)
such that  $\bK(\Fq)=G(O)$, $\bI(\Fq)=I$; and a canonical ind-group scheme
$\bG$ with $\bG(\Fq)=G(F)$. We also have  ind-varieties
$\Fl=\bG/\bI$ and $\Gr=\bG/\bK$.
 More precisely, $\Fl$, $\Gr$ are direct limits of
projective varieties with transition maps being closed imbeddings, and
$\Fl(\Fq)=G(F)/I$, $\Gr(\Fq)=G(F)/G(O)$.

The orbits of $\bI$ on $\Fl$, $\Gr$ are finite dimensional,
and isomorphic to an affine space; they are sometimes
called Schubert cells. As before,
$W_f$ is the Weyl group of $G$, and $W$ is its affine Weyl group.
Then $W$ is identified with the set of Schubert cells in $\Fl$.
For $w\in W$ (respectively $w\in W/W_f$)
let  $\Fl_w$, $\Gr_w$
be the corresponding Schubert cells.

Our main character is the derived category of
 $l$-adic sheaves
on $\Fl_{\Fqbar}$, constant along the stratification by Schubert cells;
here the subscript $\ _{\Fqbar}$ denotes extension of scalars from $\Fq$
to the algebraic closure $\Fqbar$.
More precisely,  let $D=D(\Fl)$
be the full subcategory in the derived category of
$l$-adic sheaves
on  $\Fl_{\Fqbar}$ constant along Schubert cells.
Let $\P=\P(\Fl)\subset D(\Fl)$
  be the abelian category of perverse sheaves.
It is well known that $D^b(\P)\iso D(\Fl)$
(see e.g. \cite{BGS}, Lemma 4.4.6;
we will not use this fact below).

Let $D^I=D^I(\Fl)$
 be the category of  $\bI$-equivariant $l$-adic sheaves on $\Fl$,
 and $\P^I(\Fl) \subset D^I$ be the full subcategory of perverse sheaves.
Then the functor of forgetting the equivarinat structure $\P^I(\Fl)
\to \P(\Fl)$ is known to be a full imbedding (though the one from
$D^I$ to $D$ is not).
 Thus $\P^I$ is an abelian subcategory of $\P$,
which is not closed under extensions.
The convolution product, which we denote by $*$ defines functors
$D\times D^I \to D$, and $D^I\times D^I\to D^I$.

\subsection{Central sheaves}\label{propcen}
Let $\P_\Gr$
 be the category of $\bK$ equivariant perverse
sheaves on $\Gr$. It is known  (\cite{doistor}, see also \cite{KGB}
for an alternative proof and a generalization) that
for $X, Y\in \P_\Gr$ the convolution $X*Y\in \P_\Gr$.
Further, convolution endows $\P_\Gr$ with the structure of a monoidal category,
and it  naturally extends to a structure of a commutative rigid tensor
category with a fiber functor; the resulting
Tannakian category is equivalent to the category $Rep(\LG)$ of
algebraic representation of the Langlands dual group $\LG$ over $\Qlbar$
 (see \cite{Gi}, \cite{MV}, at least for an analogous statement over $\Ce$,
and also \cite{BD}).
We will identify $Rep(\LG)$ with $\P_\Gr$.

In \cite{KGB} a functor $Z:Rep(\LG)=\P_\Gr\to P^I(\Fl)$ was constructed.
It enjoys the following properties.

i) $\pi_*\circ Z\cong id$, where $\pi:\Fl\to \Gr$ is the
projection.

ii) (Exactness of convolution) For $\F\in \P_\Gr$, $\G\in \P$
we have $\G*Z(\F)\in \P$.

iii) (Compatibility with convolution and centrality)
$Z$ is a central functor from the tensor category $\P_\Gr$ to
the monoidal category $D^I$ in the sense of definition \ref{defcen}
 above.\footnote{Compatibilities (ii), (iii) of the Definition \ref{defcen}
above are not checked in \cite{KGB}, but nevertheless hold \cite{Gap}.}

iv) (Monodromy)
A unipotent automorphism $\mon$ of the functor $Z$ is given,
$\mon_{Z(\F)}\in Aut (Z(\F))$;
it is called the monodromy automorphism (for reasons explained in
\cite{KGB}). It satisfies
\begin{equation}\label{tensprop}
  \mon_{Z(\F*\F')}
=\mon_{Z(\F)} *\mon_{Z(\F')},
\end{equation}
where we identified $End(Z(\F*\F'))= End(Z(\F)*Z(\F'))$ by means of (iii).

\subsection{Serre quotient categories}\label{SerreQ}
The set of isomorphism classes of irreducible objects in $\P$ is in
 bijection with $W$. For $w\in W$ let $L_w$ be the corresponding irreducible
object; more presicely, we set $L_w= j_{w!*}(\underline{\Ql}
[\dim \Fl_w])$, where $j_w$ denotes the imbedding $\Fl_w\imbed \Fl$, and
$\underline{\Ql}$ is the constant sheaf.

Recall that a Serre subcategory in an abelian category is a strictly
full abelian
 subcategory closed under extensions and subquotients. If $A$ is an abelian
category, and $B$ is a Serre subcategory, then the {\it Serre quotient} $A/B$
is again an abelian category. If
 every object in $A$ has finite length, then $B$ is
uniquely specified by the set of (isomorphism classes of) irreducible
objects of $A$, which lie in $B$.

Let $\c\subset W$ be a  two-sided cell. Set $W_{<\c}=\cupl_{\c' < \c} \c'$;
 $W_{\leq \c}=\cupl_{\c' \leq \c} \c'$; here $\leq$ is the standard
partial order on the set of 2-sided cells (see \cite{Cells0}), and we
write $\c'< \c$ instead of $\c'\leq \c\ \&\  \c'\ne \c$.

For a subset $S\subset W$ let $\P_S$ denote the Serre subcategory of $\P$
whose set of (isomorphism classes) of irreducible objects
equals $L_w$, $w\in S$.
We abbreviate
$\P_{\leq \c}=\P_{W_{\leq \c}}$, $\P_{<\c}=\P_{W_{ <\c }}$.
We also set $P_S^I=P^I\cap P_S$ etc.

Let $ \P_\c^I$ denote the Serre quotient category $\P_{\leq \c}^I/\P_{<\c}^I $.

\end{section}
\begin{section}{Truncated convolution categories}
\subsection{Truncated convolution and action of the central sheaves}
 Let $D_{\leq \c}(\Fl)$,  $D_{< \c}(\Fl)$ be the
 full triangulated subcategories of $D(\Fl)$ consisting of complexes with
cohomology in, respectively, $\P_{\leq \c}$, $\P_{<\c}$.
From the definition of a two-sided cell it follows that
$D_{\leq \c}$, $D_{<\c}$ are stable
 under  convolution with any object of $D^I$, i.e.
$X\in D_{\leq \c}$, $Y\in D^I$ $\Rightarrow$ $X*Y\in D_{\leq \c}$,
and same for $D_{<\c}$.

In particular, for $X, Y\in \P_{\leq \c}^I$ we have $H^i(X*Y)\in \P_{\leq \c}
^I$, and the image of $H^i(X*Y)$ in $\P_\c^I$ depends canonically only on the
image of $X,Y$ in $\P_\c^I$. Thus the formula $$(X,Y) \mapsto H^i(X*Y) \mod
 \P_{<\c}^I$$ defines a bilinear
functor $\P_\c^I\times \P_\c^I\to \P_\c^I$.

Recall that for a two-cided cell $\c$ a non-negative integer $a(\c)$
is defined (see \cite{Cells0}); and we have $H^i (X*Y)\in \P_{<\c}^I$
for $i>a(\c)$, $X\in P_{\leq \c}^I$, $Y\in \P^I$.

For $X,Y\in \P_\c^I$ we define their truncated convolution
$X\bu Y\in \P_c^I$ by $X\bu Y=
H^{a(\c)} (X*Y) \mod \P^I_{<\c}$. For semisimple $X, Y$ this coincides with
the definition  in \cite{conje}.

Also, for $\F\in \P_\Gr$ the exact functor
$\G\mapsto \G*Z(\F)$ preserves $\P_{\leq \c}^I$ and $\P_{<\c}^I$,
 hence induces an exact functor from $\P_\c$ to itself (denoted again
by $\G\mapsto \G*Z(\F)$).

\subsection{Monoidal category $\A_\c$}
\begin{Prop}\label{Ac}
Let $\A_\c$ be the strictly full subcategory of $\P_\c^I$ consisting of
all (objects isomorphic to) subquotients of $L_w*Z(\F)$, $w\in \c$,
$\F\in \P_\Gr$. Then

a) Restriction of $\bu$ to $\A_\c\times \A_\c$ takes
values in $\A_\c$, and is exact in each variable.

b) It equips
$\A_\c$ with a structure of
 a monoidal category; the object  $\I=\oplusl_d L_d$, where
$d$ runs over the set of Duflo involutions in $\c$, is a unit object of
$(\A_\c, \bu)$.
\end{Prop}

\begin{Rem}\label{rig_c}
It seems possible to  check that the monoidal category $\A_\c$ is rigid, i.e.
for $X\in \A_\c$ the ``dual'' object $\Xcheck$ is defined together
with morphsisms $ev: X\bu \Xcheck \to \I$, and $\delta:\I\to \Xcheck \bu X$
satisfying the usual compatibilities (see e.g. \cite{De} 2.1.2).
Here $\Xcheck$ has the following geometric interpretation.
The category $\P^I$ has a canonical anti-involution $\iota :
\P^I \to (\P^I)^{op}$ induced (loosely speaking) by the morphism
 $i:\bG \to \bG$, $i: g\mapsto g^{-1}$. Then
\begin{equation}\label{Xcheck}
\Xcheck \cong {\mathbb V}(\iota(X))
\end{equation}
canonically, where $\mathbb V$ stands for Verdier duality. (We neither
use nor prove this fact here).

\end{Rem}

\proof of the Proposition. a)
It follows from the definitions that $\bu$ is right exact in each
variable. We will deduce that it is exact on $\A_\c$
from a result of Lusztig on asymptotic Hecke algebras. (The argument
will use this result of Lusztig and  mid-exactness of $\bu$, but not
its  right exactness; see also Remark \ref{op1} below).

The group
$K(\P_\c)$ has a natural structure of an associative algebra; the product
on $K(\A_\c)$ is denoted by $\bu$ and is defined by
$$[L_{w_1}]\bu [L_{w_2}]=[L_{w_1} \bu L_{w_2}]$$
for irreducible objects $L_{w_1}$,  $L_{w_2} \in \A_\c$.
Thus $K(\A_\c), \bu$ is the asymptotic Hecke algebra $J_\c$, cf. \cite{conje}.
The classes of simple objects $L_w\in \P_\c$ form a standard
basis of this algebra; the standard notation for elements of this basis
is $t_w=[L_w] \in J_\c=K(\P_\c)$.

The key step is the next

\begin{Lem}\label{bubu}
For $\G_1,\G_2\in \P_\Gr$ and simple objects $L_{w_1}$,  $L_{w_2} \in \A_\c$
we have
\begin{equation}\label{eqbubu}
[(Z(\G_1)*L_{w_1})\bu (L_{w_2}*Z(\G_2))]=
[Z(\G_1)*L_{w_1}]\bu [L_{w_2}*Z(\G_2)].
\end{equation}
\end{Lem}
\proof The Lemma is a consequence of \cite{Cells2}, 2.4. Let us recall
the statement of {\it loc. cit.}

To formulate it we need some notations. Set  $A=\Zet[v,v^{-1}]$;
thus the  affine Hecke
algebra $\H$ is an $A$-algebra.
Let  $C_w\in \H$, $w\in W$ be the Kazhdan Lusztig basis of $\H$.
 Let also
$\H_{\leq \c}$, $\H_{<\c}$ be $A$ submodules in $\H$ generated by
$C_w$ with $w$ running over the corresponding subset of $W$.
According to the definition of a two-sided cell the $A$-submodules
$\H_{\leq \c}$, $\H_{<\c}$ are two-sided ideals in $\H$. Thus
$\H_\c:= \H_{\leq \c}/\H_{<\c}$ is a bimodule over $\H$.
Let $S: \H_\c \to J_c\otimes _\Zet A$ be the isomorphism of $A$-modules,
 which sends
$C_w$ into $t_w=[L_w]\in J_\c\subset J_\c\otimes A$.
Let us transport the structure  of an $\H$-bimodule to  $J_\c\otimes A$
by means of $S$;
 denote the resulting action  by $h_1\otimes h_2: x \mapsto h_1*x*h_2$,
$h_1,h_2\in \H$, $x\in J_\c\otimes A$.

\begin{Fact}\label{bimod} (\cite{Cells2}, 2.4)
a) Thus defined right
(respectively, left) action of $\H$
on $J_\c\otimes A$ commutes with the canonical left (resp., right)
action of $J_\c$
on $J_\c\otimes A$. In other words,  extending
the $\bu$ product to an $A$-algebra structure on $J_\c\otimes A$,
we get
\begin{equation}\label{actcomm}
h_1* (x\bu y)*h_2 =(h_1* x)\bu (y*h_2).
\end{equation}

b) The map $\phi_\c:\H \to J_\c \otimes A$ defined by
\begin{equation}\label{phi}
\phi_\c(h) = h*(\suml t_d),
\end{equation}
(where $d$ runs over the set of Duflo involutions in $\c$) is an algebra
homomorphism. \epf
\end{Fact}

We deduce Lemma \ref{bubu} from \eqref{actcomm}. In fact we will
use only specialization of \eqref{actcomm} at $v=1$, and only for
$h$ in the center of $\H$.

 Fix a homomorphism $A\to \Zet$ sending $v$ to 1, and set
$H=\H\otimes _A \Zet\cong
\Zet[W]$. The structure of $\H$-bimodule
 on $J_\c\otimes A$ yields a structure of
$H$-bimodule
on $J_\c$, which we also denote by $h_1\otimes h_2: x \mapsto h_1*x*h_2$.
Thus \eqref{actcomm} holds for $x,y\in J_\c$ and $h_1,h_2\in H$.

We have a standard isomorphism $K(\P)=K(\P^I)=K(D^I)\cong H=\Zet[W]$
compatible with the ring structure, where the product in $K(D^I)$ is defined
by means of convolution: $[X]*[Y]=[X*Y]$.

Let $x=t_{w_1}=[L_{w_1}]$, $y=t_{w_2}=[L_{w_2}]\in K(\P_\c^I)=J_\c$;
and $h_i=[Z(\G_i)]\in K(\P^I)=H$, $i=1,2$;
then we claim that the left-hand side of
\eqref{eqbubu} equals the left-hand side
of \eqref{actcomm}, and the right-hand side of
\eqref{eqbubu} equals the right-hand side of \eqref{actcomm}.

Indeed, the statement about the right-hand sides follows from the definitions.


To check the statement about the left-hand sides we rewrite
\begin{multline*}
(Z(\G_1)*L_{w_1})\bu (L_{w_2}*Z(\G_2))=
H^a(Z(\G_1)*L_{w_1} *L_{w_2}*Z(\G_2))\mod \P_{<\c}^I
=\\
H^a(Z(\G_1)*(L_{w_1}*L_{w_2})*Z(\G_2)) \mod \P_{<\c}^I.
\end{multline*}
By \ref{propcen} (ii) above,
 convolution with $Z(\G)$ is exact for $\G\in \P_\Gr$,
thus the right hand side in the last equality
  equals $Z(\G_1)*H^a(L_{w_1}*L_{w_2})
*Z(\G_2)$, hence its class equals $h_1*(t_{w_1}\bu t_{w_2}) *h_2$. \epf

\begin{Rem}\label{op1}
It may be possible to  get an alternative proof of exactness of $\bu$ product,
 not appealing
to Lusztig's result \eqref{actcomm}, by
establishing rigidity in $\A_\c$ by direct geometric considerations as in
Remark \ref{rig_c} above. (Recall that for an abelian tensor category
rigidity implies exactness of tensor product, cf. \cite{DM}, Proposition
1.16.)
\end{Rem}

The last Lemmas together with Lemma \ref{tochno} imply
 exactness of $\bu|_{\A_\c \times \A_\c}$
in each variable, and part (a) of the Proposition.

b)  Associativity of truncated convolution follows from
associativity of convolution, and equality $H^i(X*Y) = 0 \mod \P_{<\c}^I$
for $i>a(\c)$, $X,Y\in \P_{\leq \c}$, because
$$(X \bu Y)\bu Z\cong H^{2a}(X*Y*Z) \mod \P_{<\c}^I \cong X\bu (Y\bu Z)$$
canonically. Also, the properties of the functor $Z$ imply that
$$(L_{w_1}*Z(\G_1))\bu (L_{w_2}*Z(\G_2))\cong
(L_{w_1}\bu L_{w_2})*Z(\G_1*\G_2)$$
canonically. The right hand side of the last equality lies in $\A_\c$;
also,  exactness of $\bu|_{\A_\c\times \A_\c}$ implies that if
$X_i\in \P_\c^I$ is a subquotient of $L_{w_i}*Z(\G_i)$, $i=1,2$, then
$X_1\bu X_2$ is a subquotient of $(L_{w_1}\bu L_{w_2})*Z(\G_1*\G_2)$.
Thus $\A_\c$ is stable under the $\bu$-product. Let us
 check that $\I$ is a unit object in $\A_\c$.
Lusztig proved  (see \cite{conje}, 2.9)
that
\begin{equation}\label{I2}
\I\bu \I \cong \I;
\end{equation}
\begin{equation}\label{IL}
\I\bu L_w \cong L_w,\ \ \ \ \ \ \  \  \ w\in \c.
\end{equation}

Given \eqref{I2} we have only
 to show that $X\mapsto \I \bu X$ is an auto-equivalence of
$\A_\c$ (cf definition of a unit object of a monoidal category on p. 105
of \cite{DM}). This functor is injective on morphisms, since it is exact
and kills no irreducible objects by \eqref{IL}.
 Thus to prove it is an equivalence it suffices to construct
for all $X$ an isomorphism
\begin{equation}\label{IX}
\I\bu X\cong X,
\end{equation}
because then  $X\mapsto \I \bu X$ is surjective on isomorphism classes of
objects; and also the map $Hom(X,Y)\to Hom(I\bu X, I\bu Y)$ is an
 an injective map between vector spaces of the same dimension,
hence an isomorphism.
(In fact, the isomorphism \eqref{IX} which will be constructed
satisfies the natural  compatibilities, but we
neither check nor use this fact).

We now construct \eqref{IX}.
We know it exists for semisimple $X$ by \eqref{IL}.
It follows readily that it is also true for $X= L*Z(\G)$, where $L$
is semisimple, and  $\G\in \P_\Gr$; indeed, we have
\begin{equation}\label{ILZ}
  \I\bu(L*Z(\G))\cong
(\I\bu L)*Z(\G)\cong L_w*Z(\G).
\end{equation}
Hence it suffices to check the following. Suppose that $\iota:Y
\imbed  L*Z(\G)$ is a subobject (where $L$ is semisimple).
We need to see that
 images of the two imbeddings
$$Y\overset{\iota}{\imbed}  L*Z(\G)\overset{\eqref{ILZ}}{=}
 \I\bu( L*Z(\G));$$
$$\I\bu Y \overset{\I\bu \iota}{\imbed} \I\bu( L*Z(\G))$$
coincide. Since the  functor $?\mapsto \I\bu ?$ is injective on morphisms,
it is enough to ensure coincidence of two subobjects
of $\I\bu(\I\bu( L*Z(\G)))$: image of $\I\bu \iota$ and image of
$\I\bu(\I\bu \iota)$. The latter follows from associativity of $\bu$.
\epf

\subsection{Tannakian category $\A_d$}
Let $d\in \c$ be a Duflo involution. Then by \cite{conje}, 2.9
we have
\begin{equation}\label{LdLd}
L_d\bu L_d \cong L_d.
\end{equation}

Let $\A_d\subset \A_\c\subset \P_\c^I$ be the strictly full subcategory
consisting of all (objects isomorphic to) subquotients of $L_d*Z(\F)$,
$\F\in \P_\Gr$. Let a functor $Res_d: Rep(\LG )= \P_\Gr \to \A_d$
be defined by $Res_d(\G)=L_d*Z(\G)$.

\begin{Lem}\label{previ}
a) $\A_d\subset \A_\c$ is a monoidal subcategory,  and $L_d\subset \A_d$
is a unit object.

b) $Res_d$ has a natural structure of a central functor.

c) $\mon$
 induces a tensor automorphism of $Res_d$ (to be denoted by $\mon_d$).
\end{Lem}

\proof a) The first statement follows from
$$(L_d*Z(\G_1))\bu (L_d *Z(\G_2)) \cong (L_d\bu L_d) *Z(\G_1*\G_2)
\cong L_d * Z(\G_1*\G_2)$$
and exactness of $\bu|_{\A_\c}$. In view of \eqref{LdLd}, in order
to check that $L_d$ is a unit object we have only to show that
$$L_d\bu X\cong X$$
for $X\in \A_d$ (cf the proof of Proposition \ref{Ac} (b)).
It also follows from \cite{conje}, 2.9
that $$L_{d'}\bu L_{d}=0$$
if $d'\ne d$ are different Duflo involutions in $\c$. Hence
$L_{d'}\bu (L_d*Z(\G))\cong (L_{d'}\bu L_d)*Z(\G)=0$, and by exactness
of $\bu|_{\A_\c}$ it follows that $L_{d'}*X=0$ for $X\in \A_d$.
Thus for $X\in \A_d$ we have
$$X\cong \I\bu X=\oplusl_{d'\in \c} L_{d'} \bu X=L_d\bu X.$$
This proves (a).

(b), (c) follow immediately from, respectively, properties (iii) and (iv)
stated in section \ref{propcen}. \epf

\end{section}

\begin{section}{Main result}
Recall that if $G_1\supset G_2$ are algebraic groups over a base
field $k$, and $N\in G_1(k)$ is an element which commutes with $G_2$,
then $N$ defines a tensor automorphism of the restriction functor
$Res^{G_1}_{G_2}:Rep(G_1)\to Rep(G_2)$. We denote this automorphism by $Aut_N$.

\begin{Thm}\label{HdNd}
 There exists a pair $H_d, N_d$, where $H_d\subset \LG_\Qlbar$
is an algebraic subgroup, and $N_d\in \LG(\Ql)$ is a unipotent element
commuting with $N_d$; an equivalence of tensor categories
$\Phi_d: Rep(H_d)\cong \A_d$, and an isomorphism $Res^{\LG}_{H_d}
\cong \Phi_d \circ Res_d$, which intertwines the tensor automorphisms
$Aut_{N_d}$ and $\mon_d$.

The pair $(H_d,N_d)$ is unique up to conjugacy.
\end{Thm}

\begin{Rem}
 Rigidity (duality) in $\A_d$ can most probably be interpreted geometrically,
and is given by \eqref{rig_c} above.
\end{Rem}

\proof of the Theorem  follows directly from Lemma \ref{previ}
and Proposition \ref{genlem}. \epf

Below $k$ will denote an algebraically closed field, ${\rm char}(k)=0$.
As before, for an algebraic group $H$ over $k$
 we will write $\R(H)$
for its representation ring, and set $\Rk(H)=\R(H)\otimes \kk$. Thus
$\Rk(H)=\O(H)^{Ad}$ is the ring of conjugation invariant functions on $H$.
For $s\in H(\kk)$ we will denote the corresponding character of $\R(H)$
(or of $\Rk(H)$) by $\chi_s:\R(H) \to \kk$, $\chi_s([V])= Tr(s,V)$.

Recall the bijection between  two-cided cells in $W$ and unipotent conjugacy
 classes in $\LG(k)$ constructed by Lusztig in \cite{Cells4}. For
a two-sided cell $\c$ we let $N_\c\in \LG(\Qlbar)$ denote a unipotent element
in the corresponding conjugacy class.

\begin{Thm}\label{compa}
 For $d\in \c$ the conjugacy classes of elements $N_d$ and  $N_\c$ coincide.
\end{Thm}

\proof We will need a characterization of the bijection
$\c \leftrightarrow N_\c$.

We set $J_\c^\kk= J_\c \otimes _\Zet \kk$,
$\H^\kk=\H\otimes \kk$.

For a unipotent element $N\in \LG(\kk)$ fix a homomorphism $\gamma_N:
SL(2)\to \LG$ defined over $\kk$, such that $N=\gamma (E)$,
where $E\in SL(2,\kk)$ is the standard unipotent element. Let $s^v_{SL(2)}
\in SL(2,\kk[v,v^{-1}])$ be the diagonal matrix with entries $v,v^{-1}$; and
define an element
$s^v_N\in \LG(\kk[v,v^{-1}])$ by
$$s^v_N=\gamma_N(s^v_{SL(2)}).$$

Recall the
 homomorphism $\phi_\c: \H \to J_\c \otimes \Zet[v,v^{-1}]$, see  \eqref{phi}.

We will also make use of the isomorphism (due to Bernstein)
\begin{equation}\label{centrBe}
B:  \R(\LG)[v,v^{-1}]\iso Z(\H),
\end{equation}
where $Z(\H)$ is the center of $\H$.


The following statement follows from Theorem 4.2
in \cite{Cells4}, see also the proof of  Proposition  2.12 in {\em loc. cit.}

\begin{Fact}\label{fac}
%
%
a) The homomorphism $\phi_\c$ sends the center $Z(\H)$
to $Z(J_\c)\otimes_\Zet A$ where $Z(J_\c)$ is the center of $J_\c$.

b) Let $\chi$ be a character of $Z(J_\c)$. The
By extension of scalars it defines a character $\chi^v:
J_\c\otimes A \to \kk[v,v^{-1}]$.
The character $\chitil=\chi^v \circ (\phi_\c|_{Z(\H)}):
Z(\H) \to \kk[v,v^{-1}]$ defines a semisimple conjugacy class
  $\Omega_s\subset \LG(\overline{\kk(v)})$.

Then
 $\Omega_s\owns s^v_{N_\c} \cdot s$ for some $s\in Z_{\LG}(\gamma_{N_\c})$.

\end{Fact}


Notice 
that if $N$, $N'$ are non-conjugate unipotent elements, then
  $s^v_{N} \cdot s$
is not $ \LG(\overline{\kk (v)})$-conjugate to $s^v_{N'} \cdot s'$
 for any $s\in Z_{\LG(\kk)}(N)$,
 $s'\in Z_{\LG(\kk)}(N')$.

Thus,  in view of Fact \ref{fac}, to prove Theorem \ref{compa}
 it suffices to check that setting $\kk=\Qlbar$ we get
\begin{equation}\label{stilch}
\Omega_s\owns s^v_{N_d}\cdot s_\const
\end{equation}
for some $s_\const\in \LG(\kk)$.

The key step in the proof of Theorem \ref{compa} is the next Lemma.

We keep notations of Theorem \ref{HdNd}. In particular for a Duflo involution
$d\in \c$ we have a subgroup $H_d\subset \LG_{\Qlbar}$ with an equivalence
$Rep(H_d)\cong \A_d \subset \A_\c$; thus the Grothendieck group
$\R (H_d)=K(\A_d)$ is a subalgebra in $J_\c=K(\A_\c)$.

For $V\in Rep(\LG)$   the nilpotent endomorphism $\log (\mon_d)=\log (N_d)$
 yields a filtration (the Jacobson-Morozov-Deligne
filtration, see \cite{Weil2}, 1.6, and also \cite{Ja}, 4.1)
 on  $Res_d(V) \cong \Phi_d(V)$. Let $gr_i(Res^\LG_{H_d}(V) )
\in Rep(H_d)$ denote
the $i$-th associated graded subquotient of this filtration.

\begin{Lem}\label{phiBV} We have an equality in $J_\c[v,v^{-1}]$:
\begin{equation}\label{sumgr}
\phi_\c(B([V]))= \suml_d \suml _i v^i [gr_i (Res^\LG_{H_d}(V))]
\end{equation}
(where $\phi_\c$ is as in \eqref{phi}).
\end{Lem}

\proof of the Lemma.
Let $D^I_{mix}=D^I_{mix}(\Fl)$ be the
  $\bI$-equivariant $l$-adic
 derived category of sheaves on $\Fl$ with integral Frobenius weights
(more presicely, it is a full subcategory of the $\bI$-equivariant
$l$-adic derived category, consisting of such complexes that
cohomolgy of their stalks
carry a Frobenius action with integral weights).
Let $\P^I_{mix}\subset D^I_{mix}$ be the full subcategory of
perverse sheaves. Then $D^I_{mix}$ is a monoidal category;
the central sheaves $Z(\G)$ lift to $\P_{mix}$.

We have a standard homomorphism of abelian groups
 $\tau: K(D^I_{mix})\to \H$           
satisfying
$$\tau([X*Y])= \tau(X) \cdot \tau(Y).$$
For $V\in Rep(\LG)$ we have $\tau([Z(V)])=B([V])$ (see
\cite{KGB}, 1.2.1).

Then it follows from the definitions that
$$\phi_\c(B[V])=S (  \tau (Z(V) * \oplusl_d L_d) \mod \H_{<\c})
,$$
where $S:\H_{\leq \c}/\H_{<\c}\iso J_\c \otimes \Zet[v,v^{-1}]$
is the isomorphism sending $C_w$ to $t_w$.

Any $X\in \P^I_{mix}$ carries the canonical weight filtration (see \cite{BBD});
 let $\gr_i(X)$ denote the $i$-th associated graded
subquotient of this filtration. Then $[\gr_i(X)]$ lies in the $\Zet$-span
of $v^i C_w$.

 In particular,
$$S(\tau (\gr_i(Z(V)*L_d)) \mod \H_{<\c} )\in v^i J_\c.$$
Thus to check \eqref{sumgr} it is enough to ensure that the filtration
on $\Phi_d(V)$ induced by the weight filtration on $Z(V)*L_d$ coincides with
the canonical (Deligne) filtration associated to the nilpotent endomorphism
$\log N_\c \in End( \Phi_d(V) )$.

We claim that a stronger statement holds. Namely, we claim that the canonical
filtration on $Z(V)*L_d$ associated to the
logarithm of monodormy  coincides
with the weight filtration (this implies the desired statement, as canonical
filtration associated to a nilpotent endomorphism
is compatible with passing to a Serre quotient category).
Since $Z(V)*L_d$ is obtained by nearby cycles from a pure weight zero
perverse sheaf (cf \cite{KGB}, Proposition 6),
the latter statement is a particular case of a general Theorem of Gabber
et. al. on coincidence of monodormic and weight filtrations on nearby cycles
of a pure sheaf, see \cite{Ja}, Theorem 5.1.2. \epf

\begin{Rem}
The proof of Lemma \ref{phiBV} is the only place in this note where
the theory of Weil sheaves (rather than the theory of $l$-adic sheaves
on a scheme over $\Fqbar$, which can be safely replaced by the theory of
constructible sheaves on the corresponding complex variety) is used.
\end{Rem}

\begin{Cor}\label{co}
 Let $\phi_d: Z(\H)\to \R(H_d)\otimes A \subset J_\c \otimes A$
be the homomorphism given by $z\mapsto z*t_d$. For a semisimple $s\in H_d(\kk)$
consider  the character
$\chi_s^v: \R(H_d) \to \kk[v,v^{-1}]$ obtained from $\chi_s:\R(H_d) \to \kk$
by extension of scalars.
The character $\chitil_s=\chi_s^v \circ (\phi_d):
Z(\H)=\R(\LG) \to \kk[v,v^{-1}]$ defines a semisimple conjugacy class
in  $ \LG(\overline{\kk(v)})$.

This conjugacy class contains element $s^v_{N_d} \cdot s$. \epf
\end{Cor}

\proof of Theorem \ref{compa}. Let $M$ be an irreducible $J_\c^\kk$
module on which the idempotent $t_d\in J_\c$ acts by a nonzero operator.
 The commutative
subalgebra $K(\A_d)^\kk=\Rk (H_d) \subset t_d J_\c^\kk t_d$
acts on the finite dimensional vector space $t_d M$, with
unit element $t_d$ acting by identity;
let $u\in t_d M$ be an eigen-vector.
The corresponding character
$\chi_u: \R(H_d)\to \kk$ comes from an element $s_u\in H_d(\kk)$;
by Corollary \ref{co} the center $Z(\H)$ acts on the vector
$u\in \phi_\c^*(M)$  by the character $\chi_{s_u\cdot s^v_{N_d}}$,
which shows \eqref{stilch}. \epf


\subsection{} We now restrict attention to the unique Duflo
involution in $\c\cap W^f$; we call it $d^f$.

\begin{Thm}\label{df}
a) The set of semisimple objects of $\A_{d^f}$ is
 $\{L_w\, |\, w\in W^f\cap \c \}$.

b) $H_{d^f}$ contains a maximal reductive subgroup of the centralizer $Z_{\LG}
(N_d)$.

\end{Thm}

\begin{Rem}\label{op2} A statement stronger than that of Theorem \ref{df}(b)
  follows from results of \cite{AB}
(see section \ref{konec} below for details): we show there that in fact
$H_{d^f}=Z_{\LG}(N_\c)$.
The argument in \cite{AB} does not use
the ``nonelementary'' result about asymptotic Hecke algebras
 cited in Fact \ref{Jf} below (the proof of this result in \cite{Cells4}
relies on the theory of character sheaves).
\end{Rem}

\proof of the Theorem. The set of irreducible objects of $\A_d$ consists
of those  $L_w$, which are subquotients of $Z(V)*L_d$ for some $V\in Rep(\LG)$.
Identifying $J_\c = K(\A_\c)$ we get
$[Z(V)*L_d]=B(V)*t_d|_{v=1};$ since
$B(V)*t_d=t_d*B(V)\in t_{d^f}\cdot J_\c \cdot t_{d^f}=J_\c^f$ we
see that indeed any subquotient of $Z(V)*L_d$ has the form $L_w$, $w\in \c
\cap W^f$.

To check that $L_w\in \A_{d^f}$ for all $w\in W^f\cap \c$,
 it suffices to check that
 for any proper subset
$S\subsetneq \c\cap W^f$ there exists $V$ such that $B(V)*t_d$
does not lie in the span of $t_w$, $w\in S$. This follows from the next
Lemma.

Let $\phi_\c^f:Z(\H) \to J_\c^f\otimes A$ be the homomorphism $z\mapsto
z*t_{d^f}$.


\begin{Lem}\label{phicf}
 For any proper subalgebra $J'\subset (J_\c^f)^ \kk$
we have $im( \phi_\c^f)\not \subset J'\otimes k(v)$.
\end{Lem}

We will deduce the Lemma from another result of Lusztig (see  Remark
\ref{op2} above).

The following statement follows from (the proof of) Propositions 9.4
 in \cite{Cells4}.

\begin{Fact} \label{Jf}
Let $J_\c^f\subset J_\c$ be the span of $\{ t_w\, |\, w\in W^f\cap \c\}$.

Then $J_\c^f = t_{d^f} \cdot J_\c^f \cdot t_{d^f}$ is a commutative subalgebra.

This subalgebra is canonically isomorphic to the ring $\R(Z(N_\c))\cong \R(Z(\gamma_{N_\c}))$.
 
 For an element
$s\in Z_{\LG}(\gamma_{N_\c})$ let $\chi_s:J_\c^f\cong \R(Z(\gamma_{N_\c}))\to k$ be the corresponding character,
and $\chi_s^v:J_\c^f\otimes A\to k[v,v^{-1}]$ be obtained from $\chi_s$ by extension
of scalars.
Then the character $\tilde\chi_s$ of $Rep(\LG)$ given by
 $\tilde \chi_s: z\mapsto \chi^v(\phi_c(z)\circ t_{d^f})$ corresponds to
the conjugacy class of the element $s^v_{N_\c}\cdot s$.

\end{Fact}

\proof of Lemma \ref{phicf}.
For $\zeta\in k$ consider the map $Z(\gamma_{N_\c})\to \LG$ given by $s \mapsto s\cdot s^v_{N_\c}|_{v=\zeta}$.
Let $R_\zeta: \O(\LG)^\LG\to \O(Z(\gamma_{N_\c})^{Z(\gamma_{Z(N_\c)})})$ be the induced homomorphism  between the rings
of conjugation invariant functions. It is clear that a subalgebra $J'\in
(J_\c^f)^k\cong \O(Z(N_\c))^{Z(N_\c)}$ as in the Lemma
necessarily contains the image of $R_\zeta$ for any $\zeta\in k$. On the other hand, it is not hard
to show that the images of homomorphisms $R_\zeta$, $\zeta\in k$ generate the ring $\O(Z(N_\c))^{Z(N_\c)}$. \epf



{\em Proof} of Theorem \ref{df}(b).
Let  $I_{lus}: (J_\c^f)\otimes \kk \iso
 \Rk (Z_\LG(N_\c))$
 denote the isomorphism of Fact \ref{Jf}.

By part (a) of the Theorem we have also an isomorphism
 $I_{gaitsg}: J_\c^f \iso \R(H_d)$,
and, by Corollary
\ref{co} for $s\in H_d(\kk)$ we have  $\chi^v_s \circ I_{gaitsg}\circ
 \phi_{d^f}
= \chi_{s^v_{N_{d^f}} }$.

For two non-conjugate semisimple elements $s_1$, $s_2\in Z(\gamma_{N_\c})$ the elements
$s_1s^v_{N_\c}$ and $s_2s^v_{N_\c}$ are not conjugate in $\LG$.
Thus Fact \ref{Jf} implies that
for two different characters $\chi_1,\chi_2:
J_\c^f\to \kk$ we have $\chi_1 \circ \phi_{d^f}\ne \chi_2 \circ \phi_{d^f}$.
It follows that
for $s\in H_d^{ss}$ we have $\chi_s \circ I_{lus}=\chi_s \cdot I_{gaitsg}$,
i.e. the isomorphism $ I_{gaitsg}\otimes \kk \circ I_{lus}^{-1}:
\Rk (Z(N_\c)) \to \Rk (H_d) $ coincides with the natural restriction map.
Thus statement (b) of the Theorem follows from the next Lemma.

\begin{Lem}
If a homomorphism of algebraic  groups over a field of characteristic zero
 $i:H_1\to H_2$
induces an isomorphism of varieties $H_1/\Ad H_1 \iso H_2/\Ad H_2$ then
it induces an isomorphism of maximal reductive quotients
$H_1^{red}\iso H_2^{red}$.
\end{Lem}

\proof We can assume that the base field is algebraically closed.

Since $H/\Ad H = H^{red}/\Ad H^{red}$ for an algebraic group
$H$, we
 can replace $H_1$, $H_2$ by $H_1^{red}$, $H_2^{red}$, and assume that
$H_1$, $H_2$ are reductive. It is clear that $i$ is injective (otherwise
$i$ sends a nontrivial semisimple conjugacy class in $H_1$ to identity of
$H_2$).

Let us first check that $i$ induces an
isomorphism of connected  components of identity $H_1^0\To H_2^0$.

  It is easy to see
 that for a (not necessarily connected) reductive
 group $H$ the connected component of identity in $H/\Ad H$ is described
by $(H/\Ad(H))^0=T/N_H(T)$.
 Hence $H_1$ and $H_2$ have common Cartan $T$, and
the image of the map $N_{H_i}(T)/T\to \Aut (T)$ does not depend on $i=1,2$.
It is enough to check that the image of $N_{H_i^0}(T)/T$ in $\Aut (T)$ does
  not depend on $i$ either.

So take $x\in N_{H_1}(T)$. It is easy to see that $\Ad x|_T\in im (
N_{H_i^0}(T)/T)$ iff $\Ad x|_{H_i^0}$ is an inner automorphism of $H_i^0$.
However the latter condition is equivalent to $\dim ((H_i/\Ad H_i)^x)=\dim
 ((H_i/\Ad H_i)^0)$, where $(H_i/\Ad H_i)^x$ is the connected component of
the class of $x$. (This is so because $\dim ({\frak g}^F)$, where $F$
is a generic automorphism of a reductive Lie algebra $\frak g$
in a fixed coset modulo inner automorphisms, is maximal when the coset is
trivial).  Since  $\dim ((H_1/\Ad H_1)^x)=\dim ((H_2/\Ad H_2)^x)$
we get the statement.

To finish the proof it remains to see that $i$ induces a bijection on the
set of connected components, $i:H_1/H_1^0\iso H_2/H_2^0$.
For this it is enough to prove the statement of the Lemma for finite
groups $H_1$, $H_2$. Then it is known as Jordan's Lemma, see e.g.
\cite{Jordan}, Lemma 4.6.1.
\epf

Theorem \ref{df} is proved. \epf

\begin{Cor}
The semisimple monoidal category $A_\c^f$,
 whose set of irreducible objects
is $\{L_w\,|\, w\in \c\cap W^f\}$, and the monoidal structure is provided by
truncated convolution (see \cite{conje}),
 is equivalent to the category
of representations of $Z_{\LG}^{red}(N_\c)$, the maximal reductive quotient of $Z_{\LG}(N_\c)$.
\end{Cor}

\proof Theorem \ref{df}(a) implies that $A_\c^f$ is the category of semisimple
objects in $\A_{d^f}$. The latter is identified with $Rep(H_{d^f})$
by Theorem \ref{HdNd}, thus $A_\c^f\cong Rep(H_{d^f}^{red})$.
 In view of Theorem \ref{df}(b) we have
$H_{d^f}^{red}= Z_{\LG}^{red}(N_{d^f})$,
and by Theorem \ref{compa} $N_{d^f}$ is conjugate to $N_\c$. \epf

The statement of the Corollary was conjectured in \cite{conje}, \S 3.2.
\end{section}
\begin{section}{Announcement of some further results}\label{konec}
In this final section we describe  results related to the subject of
the present note, to appear in \cite{AB}.

Set $\P^f$ denote the Serre quotient category $\P/P_{W-W^f}$
 (notations of \ref{SerreQ}); thus  the irreducible
objects of $\P^f$ are indexed by minimal length representatives of
cosets $W_f\backslash W/W_f$ (which are in bijection with dominant coweights
of $G$). $\P^f$ carries a filtration by Serre subcategories $\P_{\leq \c}^f=
\P^f\cap \P_{\leq \c}$. Above we constructed an imbedding of the
category $Rep(H_\c)$ in the Serre subquotient category $\P^f_{\leq \c}/
\P^f_{<\c}$ for a subgroup $H_\c=H_{d^f}\subset Z_{\LG}(N_\c)$, and stated
without proof (see Remark \ref{op2}) that $H_\c=Z_{\LG}(N_\c)$
(though checked that $H_\c^{red}=Z_{\LG}(N_\c)^{red}$).

Let $\N \subset \LG$ be the subvariety of unipotent elements, and
 $Coh^\LG(\N)$ be the category
of $\LG$-equivariant coherent sheaves on $\N$.
For $N\in \N$ let $Coh_{\leq N}^\LG(\N)\subset Coh^\LG(\N)$ be the subcategory
of sheaves whose support is contained (as a set, not necessarily as a scheme)
in the closure of the orbit $\LG(N)$, and let $Coh_{< N}^\LG(\N)
\subset Coh_{\leq N}^\LG(\N)$ consist of those sheaves whose restriction to
$\LG(N)$ iz zero. Then the Serre quotient
$Coh_{\leq N}^\LG(\N)/Coh_{< N}^\LG(\N)$ is identified with the category
of equivariant sheaves on the formal neighborhood of $\LG(N)$ in $\N$; it
contains $Coh^\LG(\LG(N))=Rep(Z_\LG(N))$ as a full subcategory.

Situations described in the two previous paragraphs look similar; this suggests
a relation bewteen $\P^f$ and $Coh^\LG(\N)$. Indeed, we show in \cite{AB}
that such a relation exists on the level of derived categories.
To formulate a precise result we need another notation. Let $\Lh$
be the (abstract) Cartan subalgebra in the Lie algebra $\Lg$ of $\LG$,
and let $\Lh_0$ be the scheme of finite length $|W_f|$ defined
by $\Lh_0=\Lh\times _{\Lh/{W_f}} \{ 0\}$ (thus $\Lh_0$ is the spectrum
of the cohomology ring of the flag variety $G/B$).
 In \cite{AB} we construct an equivalence of triangulated categories
\begin{equation}\label{AB1}
D^b(\P^f)\cong D^b(Coh^\LG(\N\times \Lh_0)),
\end{equation}
 where
$\LG$ acts trivially on $\Lh_0$. Equivalence \eqref{AB1} is compatible
with one in Theorem \ref{HdNd} (for $d=d^f\in W^f$) in the natural sense;
this allows, in particular, to show that $H_{d^f}=Z_{\LG}(N_\c)$.

The tautological $t$-structure on $D^b(\P^f)$ defines, by means of
\eqref{AB1}, a $t$-structure on $ D^b(Coh^\LG(\N\times {\Lh}_0))$
(and also on $ D^b(Coh^\LG(\N))$). This $t$-structure coincides with the
perverse $t$-structure on equivariant coherent sheaves, corresponding
to the middle perversity, see Example 1 at the end of \cite{izvrat}.

Along with
the description \eqref{AB1} of the category $D^b(\P^f)$, we give a simliar
description of the ``larger''
 category $D^b(\P)$. To state we need some notations. Let
$\Lg$ be the Lie algebra of $\LG$, and let $\Ntil=T^*(\LG/\LB)\to \N$,
$\Lgtil \to \Lg$ be the Spinger-Grothendieck maps.
Set $\Stt = \Lgtil\times _\Lg \Ntil$; thus $\Stt$ is a nonreduced
scheme, whose reduced part is the ``Steinberg variety of triples''
$\Stt^{red}=St=\Ntil\times_\N \Ntil$ of $\LG$.
We construct an equivalence
\begin{equation}\label{AB2}
D^b(\P)\cong D^b(Coh^\LG(\Stt)).
\end{equation}
The isomorphism of Grothendieck groups underlying \eqref{AB2}
is the well-known isomorphism between two different realization of the
group ring $\Zet[W]$, see e.g. \cite{GC}, 7.2.
(A ``graded'' version of \eqref{AB2}, in the sense of \cite{BGS}, 4.3,
 which  relates the derived category
of mixed sheaves to the derived category of $\LG\times G_m$ equivariant
coherent sheaves on $\Stt$ gives the isomorphism between
two different geometric realizations of the affine
Hecke algebra $\H$, see \cite{GC}, chapter 7; and also
{\it loc. cit.}, Introduction, p. 15.)

The natural proper morphism $\pi:\Stt \to \N \times \Lh_0$
yields a functor $\pi_*:D^b(Coh^\LG(\Stt))\to D^b(Coh^\LG(\N\times \Lh_0))$;
under the equivalences \eqref{AB1}, \eqref{AB2} this functor corresponds
to the tautological projection functor $D^b(\P)\to D^b(\P^f)$.

\end{section}


\end{document}